\documentclass[a4paper]{amsart}
\usepackage[dvips]{graphicx}
\usepackage{verbatim,enumerate}
\usepackage{amsthm}
\theoremstyle{plain}
 \newtheorem{theorem}{Theorem}[section]
 \newtheorem{fact}{Fact}[section]
 \newtheorem*{theorem*}{Theorem}
 \newtheorem{proposition}[theorem]{Proposition}
 \newtheorem{lemma}[theorem]{Lemma}
 \newtheorem{corollary}[theorem]{Corollary}
\theoremstyle{definition}
 \newtheorem{definition}[theorem]{Definition}
\theoremstyle{remark}
 \newtheorem{remark}[theorem]{Remark}
 \newtheorem{example}[theorem]{Example}
\numberwithin{equation}{section}

\newcommand{\vect}[1]{\boldsymbol{#1}}
\newcommand{\R}{\vect{R}}
\newcommand{\C}{\vect{C}}
\newcommand{\Open}{W}

\newcommand{\rank}{\operatorname{rank}}

\newcommand{\Hyp}{\operatorname{Hyp}}
\newcommand{\PSL}{\operatorname{PSL}}
\renewcommand{\Re}{\operatorname{Re}}
\renewcommand{\Im}{\operatorname{Im}}
\newcommand{\inner}[2]{\left\langle{#1},{#2}\right\rangle}
\newcommand{\eucinner}[2]{\left\langle{#1},{#2}
          \right\rangle_{\operatorname{Euc}}}
\title[Maximal Surfaces]{
   Maximal surfaces with singularities \\
   in Minkowski space
}
\author{Masaaki Umehara}
\address[Umehara]{%
   Department of Mathematics, Graduate School of Science,
   Osaka University,
   Toyonaka, Osaka 560-0043,
   Japan
}
\email{umehara@math.wani.osaka-u.ac.jp}
\author{Kotaro Yamada}
\address[Yamada]{%
   Faculty of Mathematics,
   Kyushu University,
   Higashi-ku, Fukuoka 812-8581, Japan%
}
\email{kotaro@math.kyushu-u.ac.jp}
\begin{document}
\begin{abstract}
 We shall investigate maximal surfaces in Minkowski 3-space
 with singularities.
 Although the plane is the only complete maximal surface without singular
 points, there are many other complete maximal surfaces
 with singularities and we show that they satisfy an
 Osserman-type inequality.
\end{abstract}
\maketitle
\section*{Introduction}
It is well-known that the only complete maximal (mean curvature zero)
space-like surface in the Minkowski $3$-space $L^3$
is the plane,
and it is also well-known that any maximal surface can be locally lifted to
a null holomorphic immersion into $\C^3$
(see, for example, \cite{Kobayashi-1} or \cite{McNertney}).
However, the projection of a null holomorphic immersion to $L^3$ might
not be regular.
We shall call such surfaces {\em maxfaces\/}, and show that this class
of generalized surfaces is a rich object to investigate global geometry.

This is somewhat parallel to the case of flat surfaces in
hyperbolic $3$-space,
in which the only complete non-singular examples are
the horosphere or  hyperbolic cylinders.
But if one considers flat (wave) fronts (namely, projections of Legendrian
immersions), there are many complete examples and interesting global
properties.
See \cite{KUY-2} and \cite{KRSUY} for details.

It should be remarked that Osamu Kobayashi
(\cite{Kobayashi-1,Kobayashi-2}) gave a Weierstrass-type representation
formula for maximal surfaces and investigated such surfaces with
conelike singularities.
Using the holomorphic representation, Estudillo and Romero \cite{ER}
defined a class of maximal surfaces with singularities in more general
type, and investigated criteria for such surfaces to be a plane.
Recently, Imaizumi \cite{Imaizumi-2} studied the asymptotic behavior of
maxfaces,
and Imaizumi-Kato \cite{IK} gave a classification of
maxfaces of genus zero with at most three ends.
On the other hand, 
Lopez-Lopez-Souam \cite{LLS} classified maximal
surfaces that are foliated by circles, which includes a Lorentzian
correspondence of Riemann's minimal surface.
Fern\`andez-L\'opez-Souam \cite{FLS}
investigated the moduli space of maximal graphs over the space-like plane
with a finite number of conelike singularities.

The Lorentzian Gauss map $g$ of nonsingular maximal surface is a map
into upper or lower connected component of the two-sheet hyperboloid in
$L^3$.
By the stereographic projection from $(1,0,0)$ of
the hyperboloid to the plane, the Lorentzian Gauss map $g$
can be expressed as a meromorphic function into
$\C\cup \{\infty\}\setminus \{\zeta\in\C\,;\,|\zeta|=1\}$.
The singular set of a maxface corresponds to the set $\{|g|=1\}$,
and $g$ can be extended meromorphically on the singular set.
We shall prove in Section~\ref{sec:complete}
that a complete maxface $f\colon{}M^2\to L^3$ 
satisfies the following Osserman-type inequality
(The definition of completeness is given in Section 4.)
\[
    2\deg g \geq -\chi(M^2)+\text{(number of ends)},
\]
and equality holds if and only if  all ends are properly embedded,
where
\[
     g:M^2\longrightarrow S^2=\C\cup \{\infty\}
\]
is the Lorentzian Gauss map and $\deg g$ is its degree
as a map to $S^2$.

We also give examples for which
equality is attained
(Section~\ref{sec:complete}).
Moreover, applying the results on singularities of wave fronts in
\cite{KRSUY},
we can investigate singularities of maxfaces and 
give a criterion for a given singular point to be 
locally diffeomorphic to
cuspidal edges or swallowtails in terms of the Weierstrass data
(Section~\ref{sec:singularity}).

It should be remarked that Kim and Yang \cite{KY} very recently
constructed complete higher genus  maxfaces of two catenoidal ends. 
 Recently Ishikawa and Machida 
 \cite{IM} showed that generic singular points
   of surfaces of constant Gaussian curvature in the 
 Euclidean $3$-space and generic singular points of 
 improper affine spheres in the $3$-dimensional affine space
  are both cuspidal edges and swallowtails.
On the other hand, Fujimori \cite{Fujimori} and
Lee and Yang \cite{LY} investigate space-like surfaces 
with singularities of mean curvature
  one in the de Sitter space. (See also a forth-coming paper \cite{FRUYY}.)
   In contrast to space-like maximal surfaces,
 time-like minimal surfaces are related to Lorentz surfaces and
 a partial differential equation of hyperbolic type, see
 Inoguchi-Toda \cite{IT}.

 The authors thank Wayne Rossman, Shoichi Fujimori 
 and the referee for valuable comments.

 \section{Preliminaries}\label{sec:prelim}
 The Minkowski $3$-space $L^3$ is the $3$-dimensional affine space $\R^3$
 with the inner product
 \begin{equation}\label{eq:lorentz-metric}
     \inner{~}{~}:=-(dx^0)^2+(dx^1)^2 + (dx^2)^2,
 \end{equation}
 where $(x^0,x^1,x^2)$ is the canonical coordinate system of $\R^3$.
 An immersion $f\colon{}M^2\to L^3$ of a  $2$-manifold $M^2$ into $L^3$ is
 called {\em space-like\/} if the induced metric
 \[
     ds^2:=f^*\inner{~}{~}=\inner{df}{df}
 \]
 is positive definite on $M^2$.
 Throughout this paper, we assume that $M^2$ is orientable.
 (If $M^2$ is non-orientable, we consider the double cover.)
 Then without loss of generality, we can regard $M^2$ as a Riemann
 surface and $f$ as a conformal immersion.

 The (Lorentzian) unit normal vector $\nu$ of a space-like immersion
 $f\colon{}M^2\to L^3$ is  perpendicular to
 the tangent plane, and $\inner{\nu}{\nu}=-1$ holds.
 Moreover, it can be regarded as a map
 \begin{equation}\label{eq:normal}
     \nu\colon{}M^2 \longrightarrow H^2_{\pm}=H^2_+\cup H^2_-,
 \end{equation}%
 where
 \begin{align*}
     H^2_{+}&:=
     \bigl\{ \vect{n}=(n^0,n^1,n^2)\in L^3\,|\,
	     \inner{\vect{n}}{\vect{n}}=-1, n^0>0
     \bigr\},\\
     H^2_{-}&:=
     \bigl\{ \vect{n}=(n^0,n^1,n^2)\in L^3\,|\,
	     \inner{\vect{n}}{\vect{n}}=-1, n^0<0
     \bigr\}.
 \end{align*}
 The map $\nu\colon{}M^2\to H^2_{\pm}$ is called the {\em Gauss map\/}
 of $f$.
 A space-like immersion $f\colon{}M^2\to L^3$ is called {\em maximal\/}
 if and only if the mean curvature function vanishes identically.
 The composition of the Gauss map to the stereographic
 projection $\pi:H^2_{\pm}\to \C\cup \{\infty\}$ 
 from the north pole $(1,0,0)$ is expressed by
 \begin{equation}\label{def:g}
 g:=\pi\circ \nu=-
 \frac{\partial f^0}{\partial f^1-\sqrt{-1}\partial f^2}
 \end{equation}
 which is a meromorphic function when $f=(f^0,f^1,f^2)$ is maximal.
 We also call $g$ the {\em Gauss map\/} of $f$.
 Since $\nu$ is valued on the set $H^2_{\pm}$,
 $|g|\neq 1$ holds on $M^2$.
 The original Gauss map $\nu$ of the maximal surface
 as in \eqref{eq:normal} is rewritten by
 \begin{equation}\label{eq:gauss-map}
       \nu = \frac{1}{1-|g|^2}\bigl(
		 -(1+|g|^2),
		 2 \Re g, 2 \Im g
		   \bigr).
 \end{equation}

  A holomorphic map
 \[
      F=(F^0,F^1,F^2) \colon{}M^2\longrightarrow \C^3
 \]
  of a Riemann surface $M^2$ to $\C^3$ is called a
  {\em Lorentzian null map\/} if
 \[
     \inner{dF}{dF}=-(dF^0)^2+(dF^1)^2+(dF^2)^2 =0
 \]
  holds on $M^2$,
  where we denote by $\inner{~}{~}$ the complexification
  of the Lorentzian metric \eqref{eq:lorentz-metric}.
 Let 
 $g\colon{} M^2\to \C\cup\{\infty\}$
 be a Gauss map of the conformal spacelike maximal immersion $f$, 
then the holomorphic map 
  \[
     F:=
	\frac{1}{2}
	   \int_{z_0}^{z}
		 \bigl(
		   -2 g,(1+g^2),\sqrt{-1}(1-g^2)
		 \bigr)\omega
  \]
  is a Lorentzian null map
  defined on the universal cover $\widetilde M^2$ of $M^2$,
  and $f=F+\overline F$ holds,
 where $\omega$ is a holomorphic $1$-form on $M^2$ given by
 \[
 \omega:=\partial f^1-\sqrt{-1}\partial f^2.
 \]
 Moreover
 \begin{equation}\label{1-st}
     ds^2 = (1-|g|^2)^2\,\omega\bar \omega
 \end{equation}
 holds (see \cite{Kobayashi-1}).
 Let $ds^2_{\Hyp}$ be the hyperbolic metric on
 $\C\cup\{\infty\}\setminus\{\zeta\in\C\,;\,|\zeta|=1\}$:
 \[
    ds^2_{\Hyp} = \frac{4\,d\zeta\,d\bar \zeta}{(1-|\zeta|^2)^2}.
 \]
 Then we have
 \begin{lemma}\label{lem:pull-back-2}
  The pull-back of the metric $ds^2_{\Hyp}$ by the Gauss map $g$
  satisfies
 \begin{equation}\label{eq:gauss-curv}
     K_{ds^2}\,ds^2 = g^*ds^2_{\Hyp}=\frac{4\,dg\,d\bar g}{(1-|g|^2)^2},
 \end{equation}
  where $K_{ds^2}$ is the Gaussian curvature of $ds^2$.
  In particular, the Gaussian curvature of maximal surface in $L^3$
  is non-negative.
 \end{lemma}

 \begin{remark}\label{rmk:complete}
 {\it The only complete maximal space-like immersion is a plane.}
 This classical fact is easily proved as follows:
 Without loss of generality, we may assume $M^2$ is connected and
  simply-connected.
  Moreover, we may assume  the Gauss map is valued in $H^2_-$, that is
  $|g|<1$.
  Suppose $M^2$ is biholomorphic to the unit disk $D^2$.
  Since $(1-|g|^2)^2\omega\bar\omega<\omega\bar\omega$,
  the metric $\omega\bar\omega$ is a complete flat metric on $D^2$,
  which is impossible. 
  So $M^2$ is conformally equivalent to $\C$, then $g$ is constant.
  This implies that the image of $f$ is a plane.
 \end{remark}

 \section{Maxfaces}\label{sec:maxface}
 \begin{definition}\label{def:maximal-map}
  A smooth map $f\colon{}M^2\to L^3$ of an oriented $2$-manifold $M^2$
  into $L^3$ is called a {\em maximal map\/} if there exists an open dense
  subset $\Open\subset M^2$ such that  $f|_{\Open}$ is a maximal
  immersion.
  A point $p$ where $ds^2$ degenerates is called a {\em singular point\/}
  of  $f$.
 \end{definition}
 \begin{definition}\label{def:maxface}
  Let $f\colon{}M^2\to L^3$ be a maximal map which gives a maximal
  immersion on $\Open\subset M^2$, and $p\in M^2\setminus\Open$ a
  singular point.
  Then $p$ is called an {\em admissible singular point\/} if
  \begin{enumerate}
   \item\label{item:maxface-1}
	On a neighborhood $U$ of $p$, there
	exists a $C^1$-differentiable 
	function $\beta\colon{}U\cap\Open\to\R_+$
	such that 
	the Riemannian metric $\beta\,ds^2$ on $U\cap\Open$
	extends to a $C^1$-differentiable Riemannian metric on
	$U$, and
   \item\label{item:maxface-2}
	 $df(p)\neq 0$
  \end{enumerate}
  hold.
  A maximal map $f$ is called a {\em maxface\/} if all
  singular points are admissible.
 \end{definition}
  The condition ``$df(p)\neq 0$''\, is equivalent to ``$\rank df=1$'' at
  the  singular point $p$.
 \begin{proposition}\label{prop:max-lift}
  Let $M^2$ be an oriented $2$-manifold and  $f\colon{}M^2\to L^3$
  a maxface which is a maximal immersion on 
  an open dense subset $\Open\subset M^2$.
  Then there exists a complex structure of $M^2$ which satisfies the
  following{\rm :}
  \begin{enumerate}
   \item $f|_{\Open}$ is conformal with respect to the complex structure.
   \item There exists a holomorphic Lorentzian null immersion
	 $F\colon{}\widetilde M^2\to \C^3$
	 such that $f\circ \pi =F+\overline F$, where 
	 $\pi:\widetilde M^2\to M^2$ is the
	 universal cover of $M^2$.
  \end{enumerate}
 \end{proposition}
  The holomorphic null immersion $F$ as above
  is called the {\em holomorphic lift\/} of the maxface $f$.
 \begin{proof}[Proof of Proposition~\ref{prop:max-lift}]
  Since the induced metric $ds^2=f^*\inner{~}{~}$ gives a Riemannian
  metric on $\Open$, it induces a complex structure on $\Open$.
  So it is sufficient to construct a complex coordinate on a
  neighborhood of an admissible singular point which is compatible to the
  complex structure on $\Open$.

  Let $p$ be an admissible singular point of $f$ and $U$ a sufficiently
  small neighborhood of $p$.
  By definition, there exists a function $\beta$ on $U\cap\Open$
  such that $\beta\,ds^2$ extends to a 
 $C^1$-differentiable
  Riemannian metric on a
  neighborhood $U$. 
  We assume $U$ is simply connected.
  Then there exists a positively oriented orthonormal frame field
  $\{\vect{e}_1,\vect{e}_2\}$ with respect to $\beta\,ds^2$
  which is $C^1$-differentiable on $U$.
  Using this, we can define a $C^1$-differentiable almost complex structure $J$
  on  $U$ such that
  \begin{equation}\label{eq:complex}
       J(\vect{e}_1)=\vect{e}_2,\qquad
       J(\vect{e}_2)=-\vect{e}_1.
  \end{equation}
  Since $ds^2$ is conformal to $\beta\,ds^2$ on $\Open$,
  $J$ is compatible to the complex structure on $\Open$ induced by
  $ds^2$.
  There exists a $C^1$-differentiable decomposition
  \[
	 (T^*M^2)^{\C} = (T^*M^2)^{(1,0)}\oplus
			  (T^*M^2)^{(0,1)}
  \]
  with respect to $J$.
  Since $f$ is $C^{\infty}$-differentiable, $df$ is a smooth
  $\R^3$-valued  $1$-form. %
  So we can take the $(1,0)$-part $\zeta$ of $df$  with respect to this
  decomposition.
  Then $\zeta$ is a 
  $C^1$-differentiable $\C^3$-valued $1$-form which is
  holomorphic on $\Open$ with respect to the complex structure 
  \eqref{eq:complex}. In particular  $d\zeta$ vanishes on $W$.
  Moreover, since $W$ is an open dense subset, $d\zeta=0$ holds on
  $U$. 

  As we assumed that $U$ is simply connected,
  the Poincar\'e lemma implies that 
  there exists a $C^1$-differentiable map $F_U\colon{}U\to\C^3$  
  such that $dF_U=\zeta$.

 Since the point $p$ is an admissible singularity,
  $\zeta+\bar \zeta=df(p)\neq 0$ on $M^2$.
  In particular $\zeta\ne 0$, and 
  at least one component of $\zeta=dF_U$ does not vanish at $p$.
  If we write $F_U=(F^0,F^1,F^2)$, we can choose $j=0,1,2$ such that
  $dF^j(p)\neq 0$.
  Using this $F^j$, we define a function $z=F^j\colon{}U\to\C=\R^2$.
  Then, $z$ gives a coordinate system on $M^2$ on
  a neighborhood of $p$, because $dF^j(p)\neq 0$.
  Since  $z=F^j$ is a holomorphic function $U\cap \Open$,
  it gives a complex analytic coordinate around $p$ compatible with
  respect to that of $U\cap \Open$.
 (If $k$ is another suffix such that $dF^k(p)\ne 0$,
 then $w=F^k$ gives also a local complex coordinate system compatible with
 respect to $z$. In fact, 
 $$
 \frac{dw}{dz}=\frac{dF^k}{dF^j}=\frac{\zeta^k}{\zeta^j}
 $$
 is holomorphic on $U\cap \Open$, and satisfies the Cauchy-Riemann equation
 on $U$, since $\zeta^k,\zeta^j$ are $C^1$-differentiable and
 $U\cap \Open$ is open dense in $U$.)

  Since $p$ is arbitrary fixed admissible singularity,
  the complex structure of $\Open$ extends across each 
  singular point $p$, 
  In particular, $\partial f$ is holomorphic whole on $M^2$ and there
  exists a holomorphic map
  $F\colon{}\widetilde M^2\to \C^3$ such that $dF=\partial f$.
  Since $\partial (F+\overline F)=dF=\partial f$, 
  $F+\overline F$ differs $f\circ \pi$ by a constant. 
  So we may take $F$ such that $F+\overline F=f\circ \pi$.
  Since {$f|_W$} is a maximal immersion,
  $F$ is a {Lorentzian} null holomorphic immersion on
  $\pi^{-1}(\Open)$.
  Since $\pi^{-1}(\Open)$ is open dense subset, 
  $F$ is a {Lorentzian} null map on  $\widetilde M^2$.
  Moreover, since $df(p)\ne 0$ at each admissible singular
  point, we have %
 \[
    {dF}(q)=\partial (f\circ \pi) (q)=\partial f(p)\ne 0
  \qquad (q\in \pi^{-1}(p))
 \]
  which implies that $F$ is an immersion whole on $\widetilde M^2$.
 \end{proof}

 Conversely, a Lorentzian null immersion $F\colon{}M^2\to \C^3$
 gives a maxface $f=F+\overline F$, if it defines a maximal immersion
 on an open dense subset.
 More precisely, we have:
 \begin{proposition}\label{prop:max-proj}
  Let $M^2$ be a Riemann surface and $F\colon{}M^2\to\C^3$ a holomorphic
  Lorentzian null immersion.
  Assume
 \begin{equation}\label{eq:degenerate}
     -|dF^0|^2+|dF^1|^2+|dF^2|^2
 \end{equation}
  does not vanish identically. Then $f=F+\overline F$ is a maxface.
  The set of singularities of $f$ is points where \eqref{eq:degenerate}
  vanishes.
 \end{proposition}
 \begin{proof}
  If \eqref{eq:degenerate} does not vanish identically,
  the set
  \[
     \Open=\{-|dF^0|^2+|dF^1|^2+|dF^2|^2\neq 0\}
  \]
  is open dense in $M^2$.
  Since $F$ is Lorentzian null,
  \begin{align*}
      -|dF^0|^2+|dF^1|^2+|dF^2|^2 &=
      -|(dF^0)^2|+|dF^1|^2+|dF^2|^2 \\
     &= -|(dF^1)^2+(dF^2)^2|+|dF^1|^2+|dF^2|^2 \\
     &\geq
      -(|dF^1|^2+|dF^2|^2)+|dF^1|^2+|dF^2|^2 =0.
  \end{align*}
  Then  it holds that
  \[
       -|dF^0|^2+|dF^1|^2+|dF^2|^2 > 0 \qquad \text{on $\Open$}.
  \]
 In particular,
 $f=F+\overline F$ determines
  a conformal maximal immersion of $\Open$ into $L^3$ with induced
  metric
  \[
     ds^2 = 2\bigl(-|dF^0|^2+|dF^1|^2+|dF^3|^2).
  \]

  On the other hand, since $F$ is an immersion, $dF\neq 0$.
  Then
  \[
	 |dF^0|^2+|dF^1|^2+|dF^2|^2 > 0.
  \]
  Hence if we set
  \[
      \beta :=
	 \frac{\hphantom{-}|dF^0|^2+|dF^1|^2+|dF^2|^2}{
		   -|dF^0|^2+|dF^1|^2+|dF^2|^2}
  \]
  on $\Open$,
  $\beta$ is a positive function on $W$ such that 
  \[
     \beta\,ds^2 = 2\bigl(|dF^0|^2+|dF^1|^2+|dF^2|^2)
  \]
  can be extended to a Riemannian metric on $M^2$.
  This completes the proof.
 \end{proof}
 \begin{remark}
  Even if $F$ is a holomorphic Lorentzian null immersion,
  $f=F+\overline F$ might not be a maxface.
  In fact, for the Lorentzian null immersion
  \[
      F=(z,z,0)\colon{}\C\longrightarrow \C^3,
  \]
  $f=F+\overline F$ degenerates on whole $\C$.
 \end{remark}
 For maximal surfaces, an analogue of the Weierstrass representation
 formula is known (see \cite{Kobayashi-1}).
 Summing up, we have:
 \begin{theorem}
 [Weierstrass-type representation for maxfaces]\label{thm:maxface}
  Let $M^2$ be a Riemann surface and $f\colon{}M^2\to L^3$ a
  maxface.
  Then there exists a meromorphic function $g$ and a holomorphic $1$-form
  $\omega$ on $M^2$ such that
  \begin{equation}\label{eq:weier-maxface}
      f = \Re\int_{z_0}^{z}
		 \bigl(
		   -2 g,(1+g^2),\sqrt{-1}(1-g^2)
		 \bigr)\omega,
  \end{equation}
  where $z_0\in M^2$ is a base point.
  Conversely, let $g$ and $\omega$  be a meromorphic function and
  a holomorphic $1$-form on $M^2$ such that
 \begin{equation}\label{eq:non-deg-g}
    (1+|g|^2)^2\,\omega\bar\omega
 \end{equation}
  is a Riemannian metric on $M^2$ and $(1-|g|^2)^2$ does not vanish
  identically.
  Suppose
 \begin{equation}\label{eq:real-period-maxface}
      \Re\oint_{\gamma}
		 \bigl(
		   -2 g,(1+g^2),\sqrt{-1}(1-g^2)
		 \bigr)\omega=0
 \end{equation}
  for all loops $\gamma$ on $M^2$.
  Then  \eqref{eq:weier-maxface} defines a maxface
  $f\colon{}M^2\to L^3$.
  The set of singular points of $f$ is given by
  $\{p\in M^2\,;\,|g(p)|=1\}$. 
 \end{theorem}
 \begin{definition}
  We set
 \begin{equation}\label{eq:dsigma}
     d\sigma^2 :=
     (1+|g|^2)^2|\omega|^2
     =2( |dF^0|^2+|dF^1|^2+|dF^2|^2 ),
 \end{equation}
  and call it
  the {\em lift-metric}\/ of the maxface $f$,
  where $F=(F^0,F^1,F^2)$ is the holomorphic lift.
 \end{definition}
 The metric $\frac{1}{2}d\sigma^2$ is nothing but the pull-back
 of the canonical Hermitian metric on $\C^3$ by  the holomorphic lift
 $F$ .
 We call a pair $(g,\omega)$ in Theorem~\ref{thm:maxface}
 the {\it Weierstrass data\/} of the maxface $f$.
 As seen in \eqref{eq:gauss-map}, $g$ is the Gauss map
 on regular points of $f$.
 We also call $g\colon{}M^2\to\C\cup\{\infty\}$ the {\em Gauss map\/}
 of the maxface $f$.

 Denote by $K_{d\sigma^2}$ the Gaussian curvature of the lift-metric
 $d\sigma^2$.
 Then, by \eqref{eq:dsigma}, we have
 \begin{equation}
  \label{eq:curv-normal}
    (-K_{d\sigma^2})\,d\sigma^2 =
	    \frac{4 dg\,d\bar g}{(1+|g|^2)^2}.
 \end{equation}
 The right-hand side is the pull-back of the Fubini-Study metric of
 $P^1(\C)$ by the Gauss map $g\colon{}M^2\to \C\cup\{\infty\}=P^1(\C)$.

 \begin{remark}
  In \cite{ER}, Estudillo and Romero defined a notion of 
  {\em generalized maximal surfaces\/} as follows:
  Let $M^2$ be a Riemann surface and $f\colon{}M^2\to L^3$ a 
  differentiable map.
  Then $f$ is called a generalized maximal surface if 
  (1)  $\varphi:=\partial f/\partial z$ is holomorphic,
  (2)  $-(\varphi^0)^2+(\varphi^1)^2+(\varphi^2)^2=0$, and 
  (3)  $-|\varphi^0|^2+|\varphi^1|^2+|\varphi^2|^2$ is not identically
  zero.
  Singular points of such a surface is either 
  (A) an isolated zero of $\varphi$ (a ``branch point'') or (B) a point
  where $|g|=1$.
  Propositions \ref{prop:max-lift} and \ref{prop:max-proj} implies that
  a maxface in our sense is a generalized maximal surface without 
  singular points of type (A).
 \end{remark}

 \section{Singularities of maxfaces}
 \label{sec:singularity}

 In the previous section, we defined maxfaces as surfaces with
 singularities.
 So it is quite natural to investigate which kind of singularities appear
 on maxfaces.
 We note that $\{(x,y,z)\in \R^3\,;\, x^2=y^3\}$ is the cuspidal edge,
 and $\{(x,y,z)\in \R^3\,;\, x=3u^4+u^2v,y=4u^3+2uv,z=v\}$
 is the swallowtail.
 We shall prove the following:
 \begin{theorem}\label{thm:sing}
  Let $U$ be a domain of the complex plane $(\C,z)$ and
  $f\colon{}U \to L^3$ a maxface with the Weierstrass data
  $(g,\omega=\hat \omega\, dz)$,
  where $\hat\omega$ is a holomorphic function on $U$.
  Then
  \begin{enumerate}
   \item\label{item:sing-1}
	 A point $p\in U$ is a singular
	 point if and only if $|g(p)|=1$.
   \item\label{item:sing-2}
	 The image of $f$ around a singular point $p$ is locally
	 diffeomorphic to a cuspidal edge if and only if
	 \[
	     \Re\left(\frac{g'}{g^2\hat\omega}\right)\neq 0
	     \qquad\text{and}\qquad
	     \Im\left(\frac{g'}{g^2\hat\omega}\right)\neq 0
	 \]
	 hold at $p$, where $\omega(z)=\hat\omega(z)dz$
and $'=d/dz$.
   \item\label{item:sing-3}
	 The image of $f$ around a  singular point $p$ is locally
	 diffeomorphic to a swallowtail if and only if
	 \begin{align*}
	    \frac{g'}{g^2\hat\omega}\in\R\setminus\{0\}\qquad
	    \text{and}\qquad
	    \Re\left\{
		\frac{g}{g'}\left(
		   \frac{g'}{g^2\hat\omega}
		\right)'\right\}\neq 0
	 \end{align*}
	 hold at $p$.
  \end{enumerate}
 \end{theorem}

 In \cite{KRSUY}, a criterion for a singular point on a wave front in
 $\R^3$ to be a cuspidal edge or a swallowtail is given.
 We shall recall it and prove the theorem as an
 application of it:
 We identify the unit cotangent bundle of the Euclidean $3$-space
 $\R^3$ with 
 $\R^3\times S^2=\{(x,\vect{n})\,;\,x\in\R^3,\vect{n}\in S^2\}$,
 then
 \[
    \xi: =n_1 dx^1+n_2 dx^2+n_3 dx^3
    \qquad\bigl(x=(x^1,x^2,x^3),~\vect{n}=(n_1,n_2,n_3)\bigr)
 \]
 gives a contact form and a map
 \[
      L=(f_L,\vect{n}_L):
       U(\subset \R^2)\longrightarrow \R^3\times S^2
 \]
 is called a {\it Legendrian\/} if the pull-back of the contact form
 $\xi$ vanishes, that is $(f_L)_u$ and $(f_L)_v$ are both
 perpendicular to $\vect{n}_L$, where $z=u+\sqrt{-1}v$.
 If $L$ is a Legendrian immersion,
 the projection $f_L$ of $L$ into $\R^3$ is called a ({\em wave})
 {\em front}.

 Now let $L=(f,\vect{n}):U\to\R^3\times S^2$ be a Legendrian immersion.
 A point $p\in U$ where $f$ is not an immersion is called a
 {\em singular point\/} of the front $f$.

 By definition, there exists a smooth function $\lambda$ on $U$
 such that
 \begin{equation}\label{eq:lambda}
    f_u \times f_v =\lambda\, \vect{n}
 \end{equation}
 where $\times$ is the Euclidean vector product of $\R^3$.
 A singular point $p\in U$ is called {\it non-degenerate\/} if
 $d\lambda$ does not vanish at $p$.
 We assume  $p$ is a non-degenerate singular point.
 Then there exists a regular curve around the point $p$
 \[
    \gamma=\gamma(t):(-\varepsilon,\varepsilon)\longrightarrow U
 \]
 (called the {\em  singular curve}) such that $\gamma(0)=p$ and the image
 of $\gamma$ coincides with the set of singular points of $f$ around $p$.
 The tangential direction of $\gamma(t)$ is called
 the {\it singular direction}.
 On the other hand, a non-zero vector  $\eta\in TU$ such that
 $df(\eta)=0$ is called the {\it null direction}.
 For each point $\gamma(t)$, the 
 null direction $\eta(t)$
 determined uniquely up to  scalar multiplications.
 We recall the following
 \begin{proposition}[\cite{KRSUY}]\label{prop:KRSUY}
 Let  $p=\gamma(0) \in U$ be a non-degenerate singular point of a front
 $f:U\to \R^3$.
 \begin{enumerate}
  \item The germ of the image of the front at
	$p$ is locally diffeomorphic to a cuspidal edge 
	if and only if  $\eta(0)$
	is not proportional to $\dot\gamma(0)$,
	where $\dot\gamma=d\gamma/dt$.
  \item The germ of the image of the front at
	$p$ is locally diffeomorphic to  a swallowtail
	if and only if $\eta(0)$ is proportional to $\dot\gamma(0)$
	and
	\[
	   \left.\frac{d}{dt}\det\bigl(\dot\gamma(t),\eta(t)\bigr)
	   \right|_{t=0}\ne 0.
	\]
 \end{enumerate}
 \end{proposition}
 Now, we identify the Minkowski space $L^3$
 with the affine space $\R^3$, and denote by $\eucinner{~}{~}$
 the Euclidean metric of $\R^3$.
 To prove Theorem \ref{thm:sing}, we prepare the following:
 \begin{lemma}\label{lem:front}
  Let $f\colon{}M^2\to L^3\simeq\R^3$ be a maxface
  with Weierstrass data  $(g,\omega)$.
  Then $f$ is a projection of a Legendrian map
  $L\colon{}M^2\to \R^3\times S^2$.
Moreover, 
    $f$ is a front on a neighborhood of $p$,  and $p$ is a non-degenerate
  singular point if and only if
  \begin{equation}\label{eq:cond-front}
    \Re\left.\left(\frac{g'}{g^2\hat\omega}\right)\right|_p\neq 0,
  \end{equation}
where $\omega=\hat\omega dz$.
 \end{lemma}
 \begin{proof}
  Let $z=u+\sqrt{-1}v$ be a complex coordinate of $M^2$ around $p$
  and write $\omega=\hat\omega\,dz$,
  where $\hat\omega$ is a holomorphic function in $z$.
  Then \eqref{eq:weier-maxface} implies 
  \[
     f_z = \frac{1}{2}\bigl(
	     -2 g,1+g^2,\sqrt{-1}(1-g^2)
	   \bigr)\hat\omega,\quad
     f_{\bar z} = \frac{1}{2}\bigl( 
	    -2 {\bar g},1+{\bar g}^2,-\sqrt{-1}(1-{\bar g}^2)
	   \bigr)\overline{\hat\omega}.
  \]
  Thus, we have
  \[
    f_u\times f_v =
     -2\sqrt{-1}f_z\times f_{\bar z}
      =
      (|g|^2-1)|\hat\omega|^2(1+g\bar g,2\Re g,2\Im g),
  \]
  where $\times$ is the Euclidean vector product of $\R^3$.
  Let
  \begin{equation}\label{eq:euclidean-normal}
   \vect{n}:=\frac{1}{\sqrt{(1+|g|^2)^2+4|g|^2}}
       ( 1+ |g|^2 , 2 \Re g, 2\Im g).
  \end{equation}
Then $\vect{n}$ is the Euclidean unit normal vector of $f$,
  that is $\eucinner{df(X)}{\vect{n}}=0$ for all $X\in TM^2$,
    where $\eucinner{~}{~}$ is the Euclidean inner product.

  From now on, we assume $|g(p)|=1$, and hence $\omega(p)\neq 0$.
  At the singular point $p$, we have
  \begin{align*}
   df &= \frac{1}{2}\bigl(-2g,(1+g^2),\sqrt{-1}(1-g^2)\bigr)
         \hat\omega\,dz
       + \frac{1}{2}\bigl(
	  -2{\bar g},(1+{\bar g}^2),-\sqrt{-1}(1-{\bar g}^2)\bigr)
	\overline{\hat\omega}\,d\bar z\\
      &=\frac{1}{2}\left(
       -2,\frac{1}{g}+g,\sqrt{-1}\left(\frac{1}{g}-g\right)\right)g\omega+
      \frac{1}{2}\left(
       -2,\frac{1}{\bar g}+\bar g,
       -\sqrt{-1}\left(\frac{1}{\bar g}-\bar g\right)\right)\bar g\bar
	\omega\\
      &=(-1,\Re g,\Im g)(g\omega+\bar g\bar \omega)
  \end{align*}
  because $\bar {g}(p)=1/g(p)$.
In particular, 
  \begin{equation}\label{eq:null}
   \eta = \frac{\sqrt{-1}}{g\hat\omega},
  \end{equation}
gives the null-direction at $p$, where
we identify $T_pM^2$ with $\R^2$ and $\C$ as
  \begin{equation}\label{eq:equiv}
   \zeta=a+\sqrt{-1}b \in \C
   ~\leftrightarrow~ (a,b)\in\R^2
   ~\leftrightarrow~ a\frac{\partial}{\partial u}+
		   b\frac{\partial}{\partial  v}
   ~\leftrightarrow~ \zeta\frac{\partial}{\partial z}+
		   \bar \zeta\frac{\partial}{\partial \bar z}.
  \end{equation}

On the other hand, we have
  \[
     d\vect{n}(p) = 
     \frac{\sqrt{-1}}{2\sqrt{2}}\left(
		       \frac{dg}{g}-\frac{d\bar g}{\bar g}
		    \right) \bigl(0,\Im g,-\Re g\bigr).
  \]
  If $dg(p)=0$ then $(f,{\vect{n}})$ is not an immersion
  at $p$ because $d\vect{n}(p)=0$.
  So we may assume $dg(p)\ne 0$.
  Then the null direction of $d\vect{n}$ at $p$ is proportional
  to
  \begin{equation}\label{eq:null-n}
   \mu = \overline{\left(\frac{g'}{g}\right)}\qquad 
	\left(' = \frac{d}{dz}\right)
  \end{equation}
  under the identification with \eqref{eq:equiv}.
  On the other hand, $f$ is a front on a neighborhood at $p$ if and only 
    if $\eucinner{df}{df}+\eucinner{d\vect{n}}{d\vect{n}}$ 
    is positive definite,
    that is
    $\eta$ in \eqref{eq:null} and $\mu$ in \eqref{eq:null-n} 
  are linearly independent, or equivalently
  \[
     0\neq \det(\mu,\eta)
      =\Im(\bar\mu\eta)=\Im\frac{g'}{g}\frac{\sqrt{-1}}{g\hat\omega}
  \]
  Then we have the first part of the conclusion.
  On the other hand, the function $\lambda$ as 
 in \eqref{eq:lambda} is calculated
  as
  \begin{equation}\label{eq:lambda-maxface}
    \lambda = \eucinner{f_u\times f_v}{\vect{n}}
	    = (|g|^2-1)|\hat\omega|^2
	      \sqrt{(1+|g|^2)^2+4|g|^2}.
  \end{equation}
 Since $d\sigma^2$ in \eqref{eq:dsigma} is a Riemannian
  metric, $g$ must have pole at $p$ if $\omega(p)=0$
  for $p\in M^2$. 
 In this case, 
 $f$ is an immersion at $p$ since $|g(p)|\ne 1$.
  Hence it is sufficient to consider the case $\omega(p)\neq 0$.
  At a singular point $p$, we have
  \[
     d\lambda(p) = 2\sqrt{2}|\hat\omega|^2\left(\bar g\, dg+g\,d\bar g\right)
		  =2\sqrt{2}|\hat\omega|^2
		   \left(\frac{dg}{g}+\frac{d\bar g}{\bar g}\right)
  \]
  because $|g(p)|=1$.
  Hence $d\lambda(p)\neq 0$ if and only if $dg(p)\neq 0$.
 If \eqref{eq:cond-front} holds at $p$, $p$ is non-degenerate
  because $dg(p)\neq 0$.
  \end{proof}
 \begin{remark}\label{rem:degenerate}
  At a singular point $p$ such that 
  $g'/(g^2\hat\omega)(p)=0$, $f$ is not a front.
  For example, for a maxface $f$ defined by the Weierstrass data
  $g=e^z$ and $\omega=\sqrt{-1}dz$,
  the pair $(f(z),\vect{n}(z))$ is not an immersion 
  into $\R^3\times S^2$ at $z=0$.%
 \end{remark}
 \begin{proof}[Proof of Theorem~\ref{thm:sing}]
  We have already shown \ref{item:sing-1} in the proof
  of Lemma~\ref{lem:front}.
  Assume $\Re(g'/(g^2\hat\omega))\neq 0$ holds at a singular point $p$.
  Then $f$ is a front and $p$ is a non-degenerate singular point.
  Since the singular set of $f$ is characterized by $g\bar g=1$,
  the singular curve $\gamma(t)$ with $\gamma(0)=p$ satisfies
  $g(\gamma(t))\overline{g(\gamma(t))}=1$.
  Differentiating this, we get
  \[
      \Re\left(\frac{g'}{g}\dot\gamma\right)=0
     \qquad\left('=\frac{d}{dz},~\dot{~}=\frac{d}{dt}\right).
  \]
  This implies that $\dot\gamma$ is perpendicular 
  to $\overline{g'/g}$, that is, proportional to
  $\sqrt{-1}\,\overline{(g'/g)}$.
  Hence we can parametrize $\gamma$ as
  \[
     \dot\gamma(t) = \sqrt{-1}\overline{\left(\frac{g'}{g}\right)}
		    \bigl(\gamma(t)\bigr)
  \]
  under the identification as in \eqref{eq:equiv}.
  On the other hand, the null direction is given as \eqref{eq:null}.
  Then by Proposition~\ref{prop:KRSUY}, the germ of the image of the
  front at $p$ is locally diffeomorphic to a cuspidal edge
  if and only if $\det(\dot\gamma,\eta)\neq 0$.
  Here, 
  \[
    \det(\dot\gamma,\eta) = \Im \overline{\dot\gamma}\eta
	= -\Im\sqrt{-1}\frac{g'}{g}\frac{\sqrt{-1}}{g\hat\omega}.
  \]
  Then we have \ref{item:sing-2}.

  Next, we assume $\Im(g'/(g^2\hat\omega))=0$ holds at the singular point
  $p$.
  In this case,
  \begin{align*}
      \left.\frac{d}{dt}\right|_{t=0} \det(\dot\gamma,\eta) &=
       \Im\left(\left(\frac{g'}{g^2\hat\omega}\right)'
		 \frac{d\gamma}{dt}\right)
       = -\Re\left(\frac{g'}{g^2\hat\omega}\right)'
			\overline{\left(\frac{g'}{g}\right)}\\
      &=-\left|\frac{g'}{g}\right|^2
	\Re\left\{\frac{g}{g'}\left(\frac{g'}{g^2\hat\omega}\right)'
	   \right\}.
  \end{align*}
  Thus, the second part of Proposition~\ref{prop:KRSUY} implies 
  \ref{item:sing-3}.
 \end{proof}

 \section{Complete Maxfaces}
 \label{sec:complete}
 Firstly, we define completeness and finiteness of total curvature
 for maxfaces:
 \begin{definition}\label{def:complete}
  Let $M^2$ be a Riemann surface.
  A maxface $f\colon{}M^2 \to L^3$ is {\em complete\/}
  (resp.\ {\em of finite type}) if
  there exists a compact set $C\subset M^2$ and a symmetric $2$-tensor
  $T$ on $M^2$ such that
  $T$ vanishes on $M^2\setminus C$ and
  $ds^2+T$ is a complete metric (resp.\ a metric of finite total
  Gaussian  curvature) on $M^2$,
  where $ds^2$ is the pull-back of the Minkowski metric by $f$.
 \end{definition}

Later (Theorem\ref{thm:finite}), we shall show that
complete maxfaces are always of finite type.

 \begin{remark}
  As seen in Lemma~\ref{lem:pull-back-2},
  the Gaussian curvature of $ds^2$ is non-negative wherever
  $ds^2$ is non-degenerate.
  Then the total curvature of $ds^2 + T$ is well-defined as a real number
  or $+\infty$.
 (The total curvature of $ds^2$ itself is not well-defined
  because \eqref{eq:gauss-curv} diverges on
  the singular set $\{|g|=1\}$.
  In fact,
  the only complete maxface 
  of finite total curvature (in the sense of improper integral)
  is the plane (\cite[Theorem 5.2]{ER}).
 \end{remark}
  \begin{lemma}\label{lem:complete}
  If a maxface $f\colon{}M^2\to L^3$ is complete
  {\rm (}resp.\ of finite type{\rm )},
  then the lift-metric $d\sigma^2$ is complete
  {\rm (}resp.\ a metric of finite total absolute curvature{\rm )}
  on $M^2$.
 \end{lemma}
 \begin{proof}
 Let $(g,\omega)$ be the Weierstrass data of $f$ and
  $T$  a symmetric $2$-tensor as in Definition~\ref{def:complete}.
  Then by \eqref{1-st} and \eqref{eq:dsigma}, we have
  $ds^2+T \leq d\sigma^2$ outside the compact set $C$.
  Thus, if $ds^2+T$ is complete, so is $d\sigma^2$.

  We denote the Gaussian curvature of the metric $d\sigma^2$ by
  $K_{d\sigma^2}$.
  Then we have
 \begin{equation}\label{eq:hermgauss}
      (-K_{d\sigma^2})\, d\sigma^2 = \frac{4\,dg\,d\bar g}{(1+|g|^2)^2}
       \leq
	    \frac{4\,dg\,d\bar g}{(1-|g|^2)^2} = K_{ds^2} \,ds^2
	    \qquad\text{on $M^2\setminus C$}
 \end{equation}
  because of \eqref{eq:gauss-curv} and \eqref{eq:curv-normal}.
  Thus, if $ds^2+T$ is of finite total curvature, the total absolute
  curvature of $d\sigma^2$ is finite.
 \end{proof}

 Our definition of \lq completeness\rq \ of maxface is
 rather restrictive: %
 In fact the universal covering of complete
 maxface might not be complete since the singular set might not
 be compact on the universal cover.
 The following \lq weak completeness\rq{} seems useful in 
 some cases.
 \begin{definition}
  A maxface $f\colon{}M^2 \to L^3$ is {\em weakly complete\/}
 if the lift-metric $d\sigma^2$
 is a complete metric.
 \end{definition}

 By Lemma~\ref{lem:complete}, completeness implies weakly completeness.
 However, the converse is not true.
  For example, let
  \[
    F := \int\left(
               -2\sqrt{-1}z, 1-z^2, \sqrt{-1}(1+z^2)
             \right)
	     \frac{dz}{(z^2-1)^2}.
  \]
  Then $F$ is a Lorentzian null immersion of the universal cover of
  $\C\cup\{\infty\}\setminus\{-1,1\}$ into $\C^3$,
  and $f=F+\overline F$ gives a maxface defined on
  $\C\cup\{\infty\}\setminus\{-1,1\}$.
  Though the lift-metric
  \[
    d\sigma^2 = \frac{(|z|^2+1)^2}{|z^2-1|^2}\,dz\,d\bar z
  \]
  is complete on $\C\cup\{\infty\}\setminus\{-1,1\}$,
  the induced metric
  \[
     ds^2 = \frac{(|z|^2-1)^2}{|z^2-1|^2}\,dz\,d\bar z
  \]
  is not.
  In fact, the set of singularities (degenerate points of $ds^2$)
  is the set $\{|z|=1\}$ which accumulates at $z=\pm 1$.

 \begin{proposition}\label{prop:finiteness}
  Let $f\colon{}M^2\to L^3$ be a weakly complete maxface.
  Suppose that the lift-metric $d\sigma^2$ has finite
  absolute total curvature.
  Then the Riemann surface $M^2$ is biholomorphic to a compact Riemann
  surface $\overline{M}^2$ excluding a finite number of points
  $\{p_1,\dots,p_n\}$.
  Moreover, the Weierstrass data $(g,\omega)$ of $f$ can be extended
  meromorphically on $\overline M^2$.
 \end{proposition}

 \begin{proof}
  By our assumptions, the lift-metric $d\sigma^2$
  is a complete metric of finite absolute total curvature.
  Moreover, by \eqref{eq:hermgauss}, the Gaussian curvature of
  $d\sigma^2$ is non-positive.
  Hence by Theorem A.1 in Appendix (or 
  by Theorem~9.1 in \cite{Osserman}), $M^2$ is biholomorphic
  to $\overline M^2\setminus\{p_1,\dots,p_n\}$.

  Identifying $\C\cup\{\infty\}$ with the unit sphere $S^2$,
  the total absolute curvature of $d\sigma^2$ is nothing
  but the area of the image of the Gauss map $g\colon{}M^2\to S^2$
  counting multiplicity.
  Hence if $d\sigma^2$ is a metric of finite total curvature,
  $g$ cannot have an essential singularity at $\{p_j\}$.
  Finally, we shall prove that $p_j$ is at most a pole of $\omega$.
  If $g(p_j)\neq \infty$, there exists a neighborhood $U$ of $p_j$ in
  $\overline M^2$ such that $|g|$ is bounded on $U$.
  In this case,
  \[
      d\sigma^2 = (1+|g|^2)^2\,\omega\bar\omega \leq 
		k\,\omega\bar\omega
  \]
  holds on $U$,  where $k$ is a positive constant,
  and hence $\omega\bar \omega$ is complete at $p_j$.
  Then by Lemma~9.6 of \cite{Osserman}, $\omega$ must have a pole at
  $p_j$.

  On the other hand, if $g(p_j)=\infty$,
  $d\sigma^2\leq k(g^2\omega)\overline{(g^2\omega)}$
  holds on a neighborhood of $p_j$, where $k$ is a positive constant.
  Hence $g^2\omega$ has a pole at $p_j$.
 \end{proof}

 We call the points $p_1,\dots,p_n$ in 
 Proposition~\ref{prop:finiteness} the {\em ends\/} of the maxface $f$.
 For a weakly complete maxface of finite total absolute curvature
 with respect to $d\sigma^2$, the Gauss map $g$ is considered as 
 a holomorphic map $g\colon{}\overline M^2\to\C\cup\{\infty\}$.

 \begin{theorem}\label{thm:finite}
  If  a maxface $f\colon{}M^2\to L^3$ is 
complete, then it is of finite type.
 \end{theorem}

\begin{remark}
This assertion is essentially different from the case of minimal 
surfaces in the Euclidean 3-space.
There are many complete minimal surfaces with infinite 
total curvature like as a helicoid.
The main difference is that the Gaussian curvature of maximal
surfaces are non-negative while that of minimal surfaces are
non-positive.
\end{remark}

\begin{proof}
Since the Gaussian curvature of $f$ is nonnegative,
 Theorem 13 of Huber \cite{Hu} implies that
 $M^2$ is diffeomorphic to 
 $\overline M^2  \setminus\{p_1,\dots,p_n\}$,
  where $\overline M^2$ is a compact Riemann surface
  and $\{p_1,...,p_n\}$ is a finite subset in $\overline M^2$.
 Moreover, a modification of Theorem 15 in Huber \cite{Hu}
 yields that $M^2$ is biholomorphic to $\overline M^2
 \setminus\{p_1,\dots,p_n\}$.
 (See the introduction of Li \cite{Li} and also the  Appendix.)
 We fix an end $p_j$ arbitrary, and take a small coordinate 
 neighborhood $(U,z)$ with the origin $p_j$.
 Without loss of generality, we may assume that 
 there are no singular points on $U\setminus \{p_j\}$,
 and thus we may also assume that $|g|<1$ holds 
 on  $U\setminus \{p_j\}$ for the 
 Gauss map $g$.
 (In fact, $g$ changes to $1/g$ if we move the position of the 
 stereographic projection to the south pole.)
 By the Great Picard theorem,
 $g$ has at most pole at $z=p_j$.
 
 We now suppose $|g(p_j)|=1$ for an end $p_j$. 
 Since $g\colon{}\overline M^2\to \C\cup\{\infty\}$ is holomorphic at
 $p_j$, we can take a complex coordinate $z$  on $\overline M^2$
 such that $z(0)=p_j$ and $g(z)=a+z^k$, where $a$ is a complex number
 with $|a|=1$ and $k$ is a positive integer.
 In this coordinate, the set
 $\{z\,;\,|g(z)|^2 = (a+z^k)(\bar a+\bar z^k)=1\}$ accumulates at the
 end $z=0$.
 Thus the singular set of $f$ is non-compact, which contradicts to 
 completeness.
 Hence we have $|g(p_j)|\neq 1$.
 Then there exists a positive number $\varepsilon(<1)$
 such that $|g|^2<1-\varepsilon$
 holds on $U$.
 In this case,
 the Gaussian curvature $K_{ds^2}$ (resp.~$K_{d\sigma^2}$)
 of $ds^2$ (resp.\ $d\sigma^2$) satisfies
 \begin{equation}\label{eq:c-estimate}
    K_{ds^2}\,ds^2 = \frac{4 dg\,d\bar g}{(1-|g|^2)^2}
		   \leq \left(\frac{2}{\varepsilon}-1\right)^2
		    \frac{4\,dg\,d\bar g}{(1+|g|^2)^2}
		   =\text{const.}(-K_{d\sigma^2})\,d\sigma^2.
 \end{equation}
 Since $p_j$ is a pole of $g$,
 $d\sigma^2$ has finite total curvature
 on $U$.
 Hence $ds^2$ is of finite type at the end $p_j$.
\end{proof}

\begin{corollary}\label{cor:completeness}
A maxface  $f\colon{}M^2\to L^3$ is 
  complete if and only if
  $f$ is weakly complete  
of finite total curvature with respect to the lift-metric 
  $d\sigma^2$
  and $|g(p_j)|\neq 1$ holds for 
  each end $p_1,\dots,p_n$.
 \end{corollary}

 \begin{proof}
Let  $f\colon{}M^2\to L^3$ be a complete maxface.
  Then $f$ is of finite type by the previous theorem.
  By  Lemma~\ref{lem:complete}, $f$ is weakly complete 
  whose total absolute curvature of the lift-metric 
  is finite, and get the conclusion.
 
 Conversely, we let
  $f\colon{}\overline M^2\setminus\{p_1,\dots,p_n\}\to L^3$ be an
  weakly complete maxface whose total absolute curvature of the 
lift-metric is finite, and assume $|g(p_j)|\neq 1$ for $j=1,\dots,n$.
  Fix an end $p_j$ and assume $|g(p_j)|<1$.
  Then we can take a neighborhood $U_j$ such that $|g|^2<1-\varepsilon$
  holds on $U_j$, where $\varepsilon\in (0,1)$ is a constant.
  In this case,
  \[
     ds^2 = (1-|g|^2)^2\omega\bar\omega 
	  \geq \varepsilon^2\omega\bar\omega
	  \geq \frac{\varepsilon^2}{4}(1+|g|^2)^2\omega\bar\omega
	   =\frac{\varepsilon^2}{4}d\sigma^2
  \]
  holds on $U_j$.
  Since $d\sigma^2$ is complete at $p_j$, so is $ds^2$.
 \end{proof}

 \begin{remark}\label{imaizumi}
  In the case of  $|g(p_j)|=1$,
  the unit normal vector $\nu$ tends to a null (light-like) vector
  at the end.
 Imaizumi \cite{Imaizumi-2} 
 investigated the asymptotic behavior of such ends.
 \end{remark}

 To prove the inequality mentioned in Introduction,
 we first investigate the behavior of the holomorphic lift
 around a single end.
 \begin{proposition}\label{prop:local-osserman}
  Let $\Delta^*=\{z\in\C\,;\, 0<|z|<1\}$ and
  $f\colon{}\Delta^* \to L^3$ be a maxface
  such that an end $0$ is complete, and denote by $F$  the holomorphic
  lift of it.
  Then $dF$ has a pole at $0$ of order at least $2$.
 \end{proposition}

 \begin{proof}
  Since $f$ is complete, the lift-metric
  \[
     d\sigma^2 = 2 \left(
		      \left|\frac{dF^0}{dz}\right|^2+
		      \left|\frac{dF^1}{dz}\right|^2+
		      \left|\frac{dF^2}{dz}\right|^2
		     \right)\,
		   dz\,d\bar z
  \]
  is a complete metric at the origin.
  Then at least one of $dF^j/dz$ ($j=0,1,2$) has a pole at $z=0$.
  We assume  $dF/dz$ has a pole of order $1$ at $z=0$.
  Then  $dF/dz$ is expanded as
  \[
     \frac{dF}{dz} = \frac{1}{z}(a^0,a^1,a^2) + O(1),
  \]
  where $O(1)$ denotes the higher order terms.
  Since $f=F+\overline F$ is well-defined on a neighborhood of $z=0$,
  the residue of $dF$ at $z=0$ must be real
 (see \eqref{eq:real-period-maxface}):  $(a^0,a^1,a^2)\in\R^3$.
  On the other hand, by the nullity of $F$, we have
  \begin{equation}\label{eq:res-null}
       -(a^0)^2 + (a^1)^2 + (a^2)^2 =0.
  \end{equation}
  Here, by \eqref{def:g},
  \[
       g(0) =-\lim_{z\to 0} \frac{dF^0}{dF^1-\sqrt{-1}dF^2}
	    = -\frac{a^0}{a^1-\sqrt{-1}a^2}.
  \]
  Then by \eqref{eq:res-null},
  $|g(0)|^2 = (a^0)^2/\{(a^1)^2+(a^2)^2\} = 1$.
  This is a contradiction, because of Theorem~\ref{cor:completeness}.
  Hence $dF$ has a pole of order at least $2$.
 \end{proof}

 \begin{theorem}[Osserman-type inequality]
 \label{thm:osserman}
  Let $\overline M^2$ be a compact Riemann surface and
  $f\colon{} \overline{M}^2\setminus\{p_1,\dots,p_n\} \to  L^3$
a complete maxface.
  Then the Gauss map $g\colon{}\overline M^2\to \C\cup\{\infty\}$
  satisfies
  \[
      2\deg g \geq -\chi(M^2)+ n = -\chi(\overline M^2) + 2n,
  \]
  and equality holds if and only if all ends are properly embedded,
  that is, there exists a neighborhood $U_j$ of each end $p_j$
  such that $f|_{U_j\setminus\{p_j\}}$ is an embedding.
 \end{theorem}
 \begin{proof}
 By \eqref{eq:curv-normal}, we have
 \[
   \deg g=\frac{1}{4\pi}\int_{M^2} {(-K_{d\sigma^2})} dA_{d\sigma^2}.
 \]
  By a rigid motion in $L^3$, we may assume $g(p_j)\ne \infty$
  ($j=1,\dots,n$).
  Since
 \[
     d\sigma^2=(1+|g|^2)^2\,\omega\bar\omega
 \]
 and $\omega$ has at least pole of order $2$,
 the inequality follows from the proof of
 the original Osserman inequality for the metric
 $d\sigma^2$. (See Theorem~9.3 in \cite{Osserman}, or  \cite{Fang}).

  As we assumed $g(p_j)\neq \infty$ ($j=1,\dots,n$),
  the equality holds if and only if $\omega$ has a pole of order exactly
  $2$ at each end.
  Assume $\omega$ has pole of order $2$ at $p_j$
  and take a coordinate $z$ around $p_j$ such
  that $z(p_j)=0$.
  Without loss of generality, we  may assume $g(0)=0$.
  By a direct calculation, we have an expansion of $f(z)$ as
 \[
     f(z)  = \frac{a}{r} ( \cos\theta,\sin\theta,0)
	      + c\log r (0,0,1) + O(1)
		\qquad (z=r e^{\sqrt{-1}\theta})
 \]
  around $z=0$,
  where $a\in\R\setminus\{0\}$ and $c\in\R$ are constants
  (see \cite{KUY-0} and also \cite{Imaizumi-2}).
  If $c\neq 0$ (resp.~$c=0$),
  the end is asymptotic to the end of the Lorentzian
  catenoid as in Example~\ref{ex:catenoid} (resp.~the plane), which
  is embedded.
  Conversely, if $\omega$ has pole of order more than $2$
  at $p_j$,
 {a similar argument to that of Jorge-Meeks \cite{JM} or \cite{Schoen}
  concludes that the end is not embedded.
  (A good reference is \cite{KUY-0}.)}
 \end{proof}

 \section{Examples}

 We shall first introduce two classical examples.

\begin{figure}
 \begin{center}
 \begin{tabular}{c@{\hspace{2em}}c}
  \includegraphics[width=4cm]{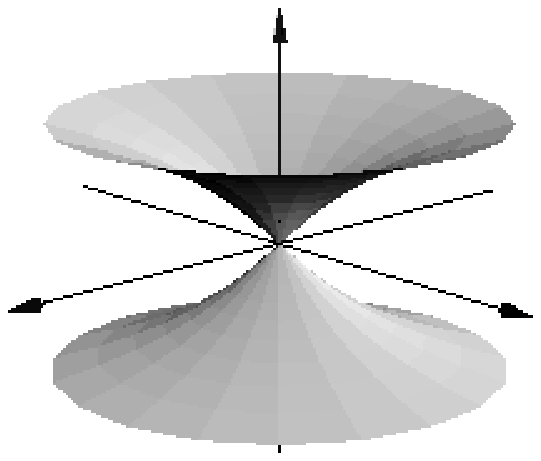} &
  \includegraphics[width=3.5cm]{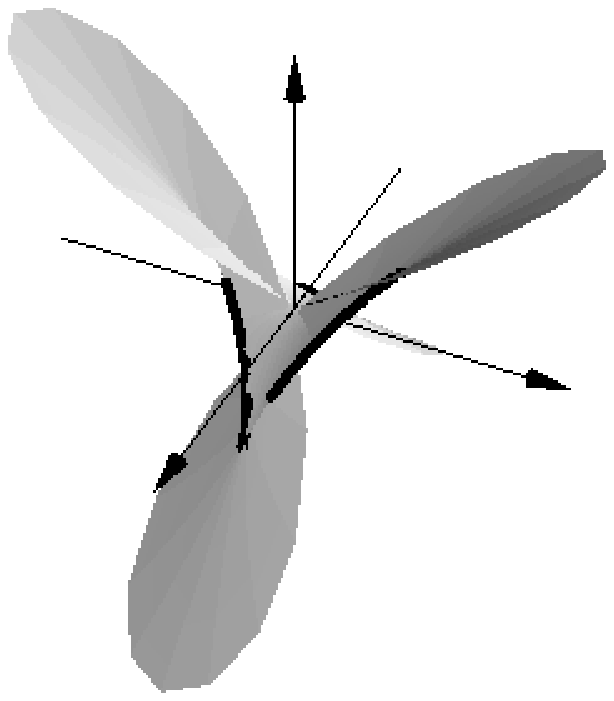} \\
  {\small the Lorentzian catenoid} &
  {\small the Lorentzian Enneper surface:}\\
   &
   {\small the singular set is shown in the black line.}
 \end{tabular}
 \end{center}
 \caption{Examples~\ref{ex:catenoid} and \ref{ex:enneper}}
 \label{fig:catenoid-enneper}
 \end{figure}

 \begin{example}[Lorentzian catenoid]
 \label{ex:catenoid}
  Rotating a curve $x^1=a\sinh(x^0/a)$ ($a\neq 0$) in the $x^0x^1$-plane
  around the $x^0$-axis,
  we have a surface of revolution
 \[
     f\colon{}S^1\times\R\ni(\theta,t)\longmapsto
	 a
	 \bigl(  t,
		 \cos\theta\sinh t,
		 \sin\theta\sinh t
	 \bigr)\in L^3.
 \]
  Then one can see that $f$ gives a maximal immersion on
  $S^1\times\R\setminus\{0\}$, hence $f$ is a maximal map
  in the sense of Definition~\ref{def:maximal-map},
  and $S^1\times\{0\}$ is the set of singularities of $f$.
  Since the induced metric is represented as
 \begin{equation}\label{eq:catenoid}
     ds^2 = a^2 \sinh^2 t(dt^2 + d\theta^2),
 \end{equation}
  $\operatorname{cosech}^2t\,ds^2=a^2(dt^2+d\theta^2)$
  extends smoothly across on the singularities.
  Hence $f$ is a maxface.
  Moreover, it can be easily seen that $f$ is
  complete. 
  The Weierstrass representation of $f$ is given
  as follows:
  Let $M^2=\C\setminus\{0\}$ and
  $g=z$, $\omega=a\,dz/z^2$.
  Then $(g,\omega)$ gives the Lorentzian catenoid
  \eqref{eq:catenoid}.
  The set of singularities is $\{|z|=1\}$ and its image by $f$ is
  the origin in  $L^3$ at which the image of $f$ is tangent
  to the light-cone
 (see Figure~\ref{fig:catenoid-enneper} left).
  Such a singularity is called a {\em conelike singularity},
  which was first investigated in \cite{Kobayashi-2}.
  See also \cite{FLS} and \cite{Imaizumi-1}.
 \end{example}
 \begin{example}[Lorentzian Enneper surface]\label{ex:enneper}
  Let $M^2=\C$ and $(g,\omega)=(z,dz)$.
  Then there exists the maxface $f\colon{}\C\to L^3$ with
  Weierstrass data $(z,dz)$.
  The set of the singularities is the unit circle $\{|z|=1\}$, and
  the points $1$, $-1$, $\sqrt{-1}$, $-\sqrt{-1}$  are swallowtails,
  and other singular points are cuspidal edges (see
  Figure~\ref{fig:catenoid-enneper}, right).
 \end{example}
 To produce further examples, we consider a relationship between
 maxfaces and minimal surfaces in the Euclidean space $\R^3$,
 and shall give a method transferring  minimal surfaces to maxfaces.

 Let $f\colon{}M^2\to L^3$ be a maxface and
 $F=(F^0,F^1,F^2)\colon{}\widetilde M^2\to \C^3$ its holomorphic
 lift,
 where $\widetilde M^2$ is the universal cover of $M^2$.
 Set
 \[
  F_0:=\bigl(F_0^1,F_0^2,F_0^3\bigr)=\bigl(F^1,F^2,\sqrt{-1}F^0\bigr).
 \]
 Since $F$ is a Lorentzian null immersion, $F_0$ is an (Euclidean) null
 immersion, that is,
 \[
       (dF^1_0)^2+
       (dF^2_0)^2+(dF^3_0)^2=0.
 \]
 Hence
 \begin{equation}\label{eq:companion}
    f_0 =  F_0+\overline{F_0}
 \end{equation}
 is a conformal minimal immersion of $\widetilde M^2$ into the Euclidean
 $3$-space $\R^3$.
 \begin{definition}
  A minimal immersion $f_0\colon{}\widetilde M^2\to \R^3$
  as in \eqref{eq:companion} is called the {\em companion\/}
  of the maxface $f$.
 \end{definition}

  The companion of the Lorentzian catenoid
  (resp.\ the Lorentzian Enneper surface)
  as in Examples~\ref{ex:catenoid} and \ref{ex:enneper}
  is the helicoid (resp.\ the Enneper surface).

 The lift-metric of a maxface $f$ as in
  \eqref{eq:dsigma} is the induced metric of the companion $f_0$,
  and the Gauss map $g_0$ of $f_0$ is represented as
 \begin{equation}\label{eq:gauss-correspondence}
   g_0 = -\sqrt{-1}g,
 \end{equation}
 where $g$ is the Gauss map of $f$.
Moreover, by Lemma~\ref{lem:complete} and by Theorem \ref{thm:finite},
 $d\sigma^2$ is a complete metric on $M^2$ with finite total curvature
 if $f$ is complete.

 By definition of $F_0$, there exists  representations
 $\rho_j\colon{}\pi_1(M^2)\to\R$ ($j=1,2,3$) such that
 \begin{equation}\label{eq:comp-period}
    F_0\circ\tau = 
       F_0 + \bigl(\sqrt{-1}\rho_1(\tau),\sqrt{-1}\rho_2(\tau),
	     \rho_3(\tau)\bigr)
     \qquad \bigl(\tau\in\pi_1(M^2)\bigr)
 \end{equation}
 holds, where $\tau$ in the left-hand side is considered as a deck
 transformation on $\widetilde M^2$.

 Conversely, we should like to construct a complete maxface
 via the complete minimal surfaces of finite total curvature:

 \begin{proposition}\label{prop:companion-inverse}
  Let $M^2$ be a Riemann surface and $\widetilde M^2$ the universal cover
  of it.  Assume a null holomorphic immersion
  $F_0\colon{}\widetilde M^2\to\C^{3}$
  satisfies the following conditions.
  \begin{enumerate}
   \item\label{item:comp-1}
	There exists representations $\rho_j$ $(j=1,2,3)$ such
	that \eqref{eq:comp-period} holds for each $\tau\in\pi_1(M^2)$.
   \item\label{item:comp-2}
	If we set
	$d F_0=(\varphi^{1}_0,\varphi^{2}_0,\varphi^{3}_0)$,
	 the function $
	    -|\varphi^{3}_0|^2+
	    |\varphi^{1}_0|^2+|\varphi^{2}_0|^2
	 $
	 does not vanish identically.
  \end{enumerate}
  Then there exists a maxface $f\colon{}M^2\to L^3$
  whose companion is $f_0=F_0+\overline{F_0}$.
  Moreover, if
  the induced metric of $f_0$ defines  a complete metric of finite
  total curvature on $M^2$,
  then  $f$ is a complete maxface
  if and only if
  \begin{equation}\label{eq:degenerate-end}
   |g(p_j)|\neq 1\qquad
   \left(
    g:=\sqrt{-1}\frac{\partial f^{3}_0}{\partial
	      f^{1}_0-\sqrt{-1}\partial f^{2}_0}\right)
  \qquad (j=1,\dots,n),
  \end{equation}
  where $M^2=\overline{M}^2\setminus\{p_1,\dots,p_n\}$
  with compact Riemann surface $\overline M^2$ and $\{p_1,\dots,p_n\}$ 
  is the set of ends.
 \end{proposition}

 \begin{proof}
  Let
 \[
     \varphi=\bigl(-\sqrt{-1}\varphi^{3}_0,\varphi^{1}_0,\varphi^{2}_0\bigr)
     \qquad\text{and}\qquad  
     F:=\displaystyle \int_{z_0}^z \varphi,
 \]
 where $z_0\in \widetilde M^2$ is a base point.
  Then by the condition \ref{item:comp-1}, $f=F+\overline F$ is
  well-defined on $M^2$.
  Moreover, by \ref{item:comp-2}, $\inner{dF}{d\overline F}$ is
  not identically $0$.
  Hence $f$ is a maxface.
  Suppose now that $f_0$ is complete and of finite total curvature.
  Then, the lift-metric
  $d\sigma^2$ is complete and  of finite total curvature.
  If $|g(p_j)|\neq 1$ for $j=1,...,n$,
  $f$ is complete by Theorem \ref{cor:completeness}.
 \end{proof}
 \begin{example}[Lorentzian Chen-Gackstatter surface]
 We set
 \[
    \overline M^2 = \{(z,w)\in\C^2\cup\{\infty,\infty\}\,;\,
	 w^2=z(z^2-a^2)\},
 \]
  where $a$ is a positive real number, and
  set
 \[
     \varphi_0:=
      \frac{1}{2}
      \left(
	\left(\frac{z}{w}-B^2\frac{w}{z}\right),
	\sqrt{-1}\left(\frac{z}{w}+B^2\frac{w}{z}\right),
	   2 B
      \right)\,dz,
 \]
  where
 \[
    B^2 =
      \left.
      \int_0^a \frac{x\,dx}{\sqrt{x(a^2-x^2)}}\right/
      \int_0^a \frac{(a^2-x^2)\,dx}{\sqrt{x(a^2-x^2)}}.
 \]
  Then
 \[
    f_0:=F_0+\overline{F_0} ,\qquad
    \left(F_0=\int_{z_0}^z \varphi_0\right)
 \]
  gives a complete minimal immersion 
 $f_0$
  of $M^2:=\overline{M}^2\setminus\{\infty\}$ into $\R^3$, %
  which is called {\em Chen-Gackstatter surface}
   (\cite{CG}).

  Since the third component of $\varphi_0$ is an exact form, 
  $F_0^3$ is well-defined on $M^2$.
  In particular, \eqref{eq:comp-period} holds for $\rho_3=0$.

 Moreover, %
  $g = Bw/z$ in \eqref{eq:degenerate-end} tends to $0$ as $z\to\infty$.
  Hence
  the corresponding maxface
  $f\colon M^2\to L^3$ given by Proposition~\ref{prop:companion-inverse}
  is complete, which is called
  the {\it Lorentzian Chen-Gackstatter surface}.
 \end{example}

 \begin{example}[Minimal surfaces which admits a Lopez-Ros deformation]
 Let $f:M^2\to \R^3$ be a complete conformal minimal immersion of
 finite total curvature. 
 Then there exists a null holomorphic lift 
 $F\colon{}\widetilde M^2\to \C^3$ such that $f=F+\overline{F}$.
 Then we have a representation 
 $\rho\colon{}\pi_1(M^2)\to \R^3$
 such that
 \begin{equation}\label{eq:comp-period-2}
    F\circ\tau = F + \sqrt{-1}\rho(\tau)
     \qquad \bigl(\tau\in\pi_1(M^2)\bigr).
 \end{equation}
 Then $f$ is called 
 {\em a minimal surface which admits the Lopez-Ros deformation\/} if
 $\rho(\pi_1(M^2))$ is contained in a $1$-dimensional subspace
 of $\R^3$.
 In this case, by a suitable rotation of the surface,
 we may assume that 
 \begin{equation}\label{eq:reducible-repr}
    \rho\bigl(\pi_1(M^2)\bigr)\subset \R\bigl(0,0,1\bigr).
 \end{equation}
 We set
 \[
    dF=\frac{1}{2}\bigl((1-g^2), \sqrt{-1}(1+g^2), 2g\bigr)\omega.
 \]
 For each non-zero real number $\lambda$, 
 replacing Weierstrass data $(g,\omega)$ by $(\lambda g,\omega/\lambda)$,
 the new minimal immersion
 \[
    f_\lambda=F_\lambda+\overline{F_\lambda},\qquad
	      dF_\lambda=\frac{1}{2}\bigl((1-\lambda^2g^2), 
		\sqrt{-1}(1+\lambda^2g^2), 2g\bigr)\omega
 \]
 also gives a conformal minimal immersion of $M^2$
 because of \eqref{eq:reducible-repr}.
 In particular, $f_{\lambda}$ is complete and of finite total curvature.
 The $1$-parameter family $\{f_\lambda\}$ is called a 
 {\em Lopez-Ros deformation of $f$ } {\cite{LR}}.  
 Then one can easily check that all of {$\sqrt{-1}F_\lambda$}
 satisfy \eqref{eq:comp-period} with $\rho_1=\rho_2=0$.  %
 Moreover, except for only finite many values of $\lambda$, the condition
 $|g(p_j)|\ne 1$ holds $(j=1,\dots,n)$, where $p_1,\dots,p_n$
 are ends of the immersion $f$.
 Thus we can construct complete maxface $\hat f_\lambda$ 
 from $\sqrt{-1}F_\lambda$ except for at most
 finitely many values of $\lambda$. 
 Remark that $f_{\lambda}$ and $f_{-\lambda}$ are congruent with 
 each other.
 The number of $\lambda>0$ such that  $\hat f_{\lambda}$ is not complete
 is not exceed the number of the ends $n$.

 Many examples of minimal surfaces which admit Lopez-Ros deformation
  are known \cite{L,KaUY,McCune}.
  So we have uncountably many examples of complete maxfaces.
 \end{example}

 \begin{example}\label{ex:jorge-meeks}
  The Jorge-Meeks surface is a complete minimal surface in $\R^3$
  with $n$ catenoidal ends.
  Such a surface is realized as an immersion
  \[
     f_0\colon{}\C\cup\{\infty\}\setminus\{1,\zeta,\dots,\zeta^{n-1}\}
      \longrightarrow \R^3\qquad
	(\zeta=e^{2\pi i/n})
  \]
  with Weierstrass data
  \[
      g_0 = z^{n-1}, \qquad \omega_0 = \frac{dz}{(z^n-1)^2}.
  \]
  It can be easily checked that,
  there exists a maxface $f$ whose minimal companion is
  the Jorge-Meeks' surface.
  Here, $g$ in \eqref{eq:degenerate-end} is $z^{n-1}$,
  $|g|=1$ holds on each end.
  Hence the maxface $f$ is not complete but weakly complete.
  As pointed out in Remark~\ref{imaizumi}, Imaizumi \cite{Imaizumi-2}
  investigated weakly complete maxface, and 
  introduced a notion of {\em simple ends\/} for 
  an end $p_j$ satisfying $|g(p_j)|=1$.
  Imaizumi and Kato \cite{IK} classified weakly complete
 maxfaces of genus zero with  at most $3$ simple ends.
 \end{example}

\appendix
\section{A consequence of Huber's theorem.}\label{app}
This appendix was prepared for the forthcoming paper
Fujimori, Rossman, Umehara, Yamada and Yang
\cite{FRUYY}, but the other authors allow to put it in this paper.
We shall show that the following assertion 
is a simple consequence of Huber's theorem
\cite[Theorem 13]{Hu}.

\begin{theorem}\label{thm:app}
 Let $(M^2,ds^2)$ be a complete Riemannian $2$-manifold and
 $K$ the curvature function. 
 Suppose 
 \begin{equation}\label{eq:finite-negative}
   \int_{M^2} (-K_-)\, dA <\infty,
 \end{equation}
 where 
 \[
   K_-:=\min(K,0).
 \]
 Then there exists a compact Riemann surface $\overline{M}^2$
 and finite points $p_1,\dots,p_n\in \overline{M}^2$ such that 
 $M^2$ is bi-holomorphic to $\overline{M}^2\setminus \{p_1,\dots,p_n\}$.
\end{theorem}

This assertion was pointed out in Li \cite{Li} without proof.
Here, we shall give a proof for
a help for readers.
To prove the assertion we use the following well-known fact
in Huber's paper

\begin{fact}[Huber {\cite[Theorem 13]{Hu}}]
\label{fact1}
 Let $(M^2,ds^2)$ be a complete Riemannian $2$-manifold
such that \eqref{eq:finite-negative} holds.
 Then $M^2$ is diffeomorphic to 
 $\overline{M}^2\setminus \{p_1,\dots,p_n\}$,
where $\overline{M}^2$ is a compact $2$-manifold.
\end{fact}

Moreover, the following assertion is known:
\begin{fact}[Blanc-Fiala \cite{BF}, Huber {\cite[Theorem 15]{Hu}}]

\label{fact2}
 Let $(M^2,ds^2)$ be a complete Riemannian $2$-manifold
 such that \eqref{eq:finite-negative} holds.
 Then $M^2$ is parabolic¡¥
\end{fact}

This assertion firstly proved by Blanc and Fiala when $M^2$ is simply
connected.
To prove our theorem, we apply the above fact 
for this simply connected case. 
\begin{proof}[Proof of the Theorem~\ref{thm:app}]
 By Fact~\ref{fact1}, $M^2$ is diffeomorphic to
 $\overline{M}^2\setminus \{p_1,\dots,p_n\}$.
 We fix an end $p_j$.
 There exists a coordinate neighborhood
 $(U_j;u,v)$ of $\overline{M}^2$ such that
 $p_j$ corresponds to the origin and the boundary $\partial U_j$
 is a simple closed $C^\infty$-regular curve.
 Then by uniformization theorem of annuli
 (see Ahlfors-Sario, Riemann surface (Princeton), I4D, II3B),
 $(U_{j}\setminus\{p_{j}\},ds^2|_{U_j\setminus\{p_j\}})$
 is conformally equivalent to 
 \[
     \Delta(r):=\{z\in \C\,;\, r<|z|<1\}
 \]
 where $r\in [0,1)$.
 Then $ds^2$ can be considered as a metric defined on $\Delta(r)$.
 Let $\lambda:\Delta(r) \to [0,1]$ be a $C^\infty$-function such that
 \begin{enumerate}
  \item $\lambda(z)=0$ when $|z|\le (r+1)/2$ 
  \item $\lambda(z)=1$ when $|z|\ge (r+2)/3$ 
 \end{enumerate}
 Then we can define a new metric $d\sigma^2$ by
 \[
   d\sigma^2=(1-\lambda)ds^2+\lambda \frac{4dz\,d\bar z}{(1+|z|^2)^2}.
 \]
 Since this metric is constant Gaussian curvature $1$
 when $|z|>1$ and can be extended at $z=\infty$,
and we obtain a complete simply connected Riemannian
 manifold
 \[
   (\{z\in \C\cup \{\infty\}\,;\,|z|>r\}, d\sigma^2)
 \]
 Moreover, the integral of the negative part of the curvature function
 of $d\sigma^2$ is finite.
 So we can apply Fact~\ref{fact2} and can conclude that
 $(\{z\in \C\cup \{\infty\}\,;\,|z|>r\}, d\sigma^2)$ is parabolic, namely
 conformally equivalent to $\C$, which implies that
 $(\{z\in \C\,;\,r<|z|<(r+1)/2\}, ds^2|_{\{z\in \C\,;\,r<|z|<(r+1)/2\}})$ 
 is conformally equivalent to
 a punctured disc. Thus we must conclude that $r=0$.
\end{proof}

 
\end{document}